\documentclass{amsart}

\usepackage{amsmath}
\usepackage{amssymb}
\usepackage{mathrsfs}
\usepackage{stmaryrd}

\usepackage{hyperref}
\usepackage{enumitem}

\usepackage{graphicx}
\usepackage[all]{xy}

\usepackage[nobysame]{amsrefs}

\newtheorem{theorem}{Theorem}[section]
\newtheorem{cor}[theorem]{Corollary}
\numberwithin{equation}{section}

\newcommand{\onto}{\twoheadrightarrow} 
\newcommand{\epsi}{\varepsilon}
\DeclareMathOperator{\Hom}{Hom} 
\newcommand{\Ext}{\mathrm{Ext}}
\newcommand{\cExt}{\mathcal{E}xt}

\DeclareMathOperator{\Spec}{Spec} 
\DeclareMathOperator{\Spf}{Spf}
\newcommand{\GL}{\mathrm{GL}} 
\newcommand{\pow}[1]{\llbracket  #1 \rrbracket} 
\newcommand{\frakm}{\mathfrak{m}}
\newcommand{\frakX}{\mathfrak{X}}
\newcommand{\conv}[1]{\mathrm{conv} \! \left\{ #1  \right\} }

\newcommand{\stack}[1]{ \mathcal{M}^\mathrm{Kss}_{#1} } 
\newcommand{\modspace}[1]{M^\mathrm{Kps}_{#1}} 

\newcommand{\DefqG}[1]{\mathrm{Def}^{\mathrm{qG}}_{#1}}
\newcommand{\Def}[1]{\mathrm{Def}_{#1}}
\newcommand{\TTqG}[2]{\mathbb{T}^{\mathrm{qG}, #1}_{#2}}
\newcommand{\cTqG}[2]{\mathscr{T}^{\mathrm{qG}, #1}_{#2}}

\newcommand\cO{\mathcal{O}}
\newcommand\cT{\mathscr{T}}
\newcommand\rH{\mathrm{H}}

\renewcommand\AA{\mathbb{A}}
\newcommand\CC{\mathbb{C}}
\newcommand\GG{\mathbb{G}}
\newcommand\PP{\mathbb{P}}
\newcommand\QQ{\mathbb{Q}}
\newcommand\RR{\mathbb{R}}
\newcommand\TT{\mathbb{T}}
\newcommand\ZZ{\mathbb{Z}}

\title[K-moduli of Fano 3-folds can have embedded points]{K-moduli of Fano 3-folds can have \\ embedded points}

\author{Andrea Petracci}
\address{Dipartimento di Matematica, Università di Bologna, Piazza di Porta San Donato 5, Bologna, 40126, Italy}
\email{a.petracci@unibo.it}

\begin{document}

\begin{abstract}
We  exhibit an example of obstructed K-polystable Fano $3$-fold $X$
such that the K-moduli stack of K-semistable Fano varieties
and the K-moduli space of K-polystable Fano varieties have an embedded point at $[X]$.
\end{abstract}

\maketitle

\section{Introduction}

We always work over an algebraically closed field of characteristic $0$, which is denoted by $\CC$. A Fano variety is a normal projective variety $X$ over $\CC$ such that its anticanonical divisor $-K_X$ is $\QQ$-Cartier and ample. In some sense, a projective variety is Fano if it has `positive curvature'.
Fano varieties occupy a prominent r\^ole in algebraic geometry from many points of view.

Moduli of Fano varieties are quite elusive, as they are highly non-separated.
Recently, K-stability \cite{tian_KE, donaldson_stability} (i.e.\ the study of the existence of K\"ahler--Einstein metrics) has been spectacularly applied to construct `reasonable' (i.e.\ separated and proper) moduli spaces of Fano varieties \cite{ABHLX, xu_minimizing, BLX, jiang_boundedness, blum_xu_uniqueness, xu_zhuang, properness_K_moduli, projectivity_K_moduli_final, lwx}.
More precisely, it is known that,
for each integer $n \geq 1$, $\QQ$-Gorenstein families of K-semistable Fano $n$-folds over $\CC$ form  an Artin stack $\stack{n}$, called the K-moduli stack, which has the following properties: $\stack{n}$ has countably many connected components and every connected component of $\stack{n}$ is of finite type over $\CC$; furthermore, the stack $\stack{n}$ has a good moduli space which is denoted by $\modspace{n}$, is called the K-moduli space, parametrises K-polystable Fano $n$-folds, and is a countable disjoint union of projective schemes over $\CC$.
We refer the reader to \cite{xu_survey} for a survey on these topics.

Since Fano varieties of dimension $2$ are unobstructed \cite{hacking_prokhorov, procams}, each connected component of the K-moduli stack $\stack{2}$ is smooth over $\CC$ and each connected component of the K-moduli space $\modspace{2}$ is a normal variety over $\CC$.

In joint work with Kaloghiros \cite{ask_petracci} we exhibited the first examples of singular points on $\stack{n}$, for each $n \geq 3$.
More precisely, we showed that, if $n \geq 3$, $\stack{n}$ and $\modspace{n}$ can be locally reducible or non-reduced (or both). These examples are constructed via toric geometry.

In \cite{petracci_murphy} we have showed that the deformation space of every isolated Gorenstein toric $3$-fold singularity appears, in a weak sense, as a singularity on $\stack{n}$, for each $n \geq 3$; we also prove that the number of local branches can be arbitrarily high.
Here we analyse the construction of \cite{petracci_murphy} in a particular example and we prove:

\begin{theorem} \label{thm}
There exists a K-polystable toric Fano $3$-fold $X$ with canonical singularities, such that
 the miniversal ring of the stack $\stack{3}$ at the point corresponding to $X$ is
	\begin{equation*}
		\CC \pow{t_1, \dots, t_8} / (t_1^2, t_1 t_2, t_3^2, t_3 t_4),
	\end{equation*}
and the completion of the structure sheaf of the scheme $\modspace{3}$ at the point corresponding to $X$ is 
\begin{equation*}
\CC \pow{u_1, \dots, u_6} / (u_1^2,  u_2^2, u_1 u_2, u_1 u_3, u_2 u_3).
\end{equation*}
\end{theorem}

An immediate consequence is:

\begin{cor}
	The stack $\stack{3}$ and the scheme $\modspace{3}$ have embedded points.
\end{cor}

It is not known whether there are any restrictions on the singularities that can appear on $\stack{n}$ or on $\modspace{n}$ for $n \geq 3$.

\subsection*{Notation and conventions}
We always work over an algebraically closed field of characteristic $0$, which is denoted by $\CC$.
We assume that the reader is familiar with toric geometry; every toric variety or toric singularity is assumed to be normal.

Throughout this article we use the following notation.

\begin{tabular}{ll}
	$P$ & the $3$-dimensional lattice polytope in Figure~\ref{fig:polytope_P} \\
	$Q$ & the polar of $P$ \\
	$\conv{\cdot}$ & convex hull of a set \\
	$X$ & the toric variety associated to the face fan of $P$ \\
	$\frakX$ & the index $1$ cover stack of $X$ \\
	$\DefqG{X}$ & the $\QQ$-Gorenstein deformation functor of $X$ \\
	$A$ & the hull of $\DefqG{X}$ \\
	$T$ & the $3$-dimensional algebraic torus acting on $X$ \\
	$M$ & the character lattice of $T$ \\
$N$ & the cocharacter lattice of $T$ \\
	$G$ & the automorphism group of $X$
\end{tabular}

\subsection*{Acknowledgements}
I am extremely grateful to Alessio Corti, Paul Hacking and Anne-Sophie Kaloghiros for numerous fruitful discussions.
I am indebted to Liana Heuberger for spotting a mistake in \S\ref{sec:mirror_symmetry} in a previous version of this article.

The titles of \S\ref{sec:polytope}, \S\ref{sec:variety}, \S\ref{sec:computation} in this article reflect the title ``Varieties, Polyhedra, Computation'' of the thematic Einstein Semester on Algebraic Geometry which was held in Berlin during the Winter Semester 2019/20: I wish to thank all participants for such a stimulating environment.

The author acknowledges partial financial support
from PRIN2020 2020KKWT53 ``Curves, Ricci flat Varieties and their Interactions'',
from ``Gruppo Nazionale per le Strutture Algebriche, Geometriche e le loro Applicazioni'' GNSAGA of INdAM,
and from the European Union -- NextGenerationEU under the National Recovery and Resilience Plan (PNRR) -- Mission 4 Education and research -- Component 2 From research to business -- Investment 1.1, Prin 2022, ``Geometry of algebraic structures: moduli, invariants, deformations'', DD N.~104, 2/2/2022, proposal code 2022BTA242 -- CUP J53D23003720006.

\section{The polytopes $P$ and $Q$} \label{sec:polytope}

\begin{figure}
	\centering
	\includegraphics[scale = 0.5]{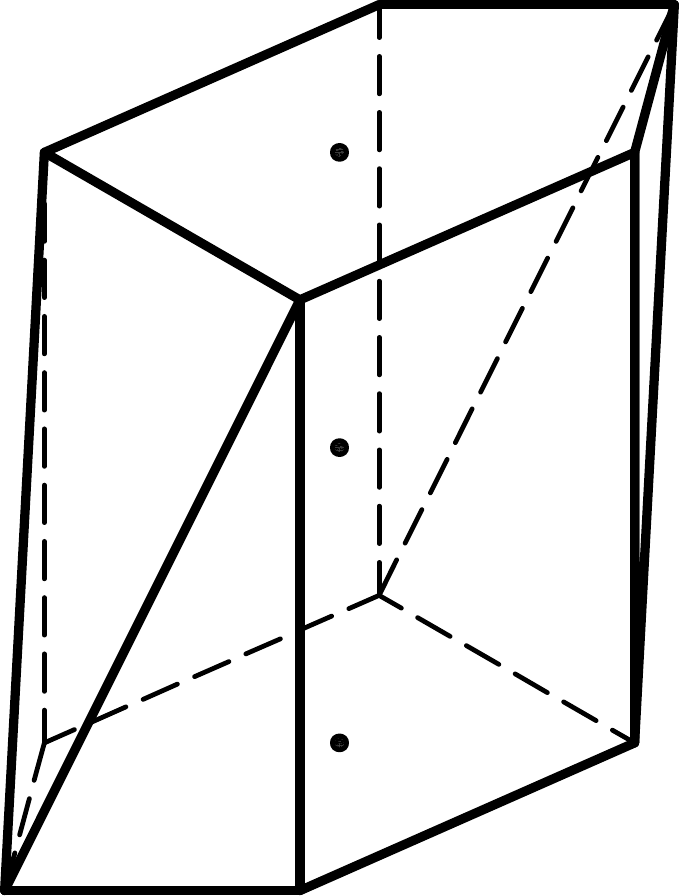}
	\caption{The polytope $P$}
	\label{fig:polytope_P}
\end{figure}

Consider the points
\begin{equation*}
	a = \! \begin{pmatrix}
		1 \\ 0 \\ 1
	\end{pmatrix} \quad
	b=  \! \begin{pmatrix}
		1 \\ 1 \\ 1
	\end{pmatrix} \quad
	c = \! \begin{pmatrix}
		0 \\ 1 \\ 1
	\end{pmatrix} \quad
	d= \! \begin{pmatrix}
		-1 \\ 0 \\ 1
	\end{pmatrix} \quad
	e= \! \begin{pmatrix}
		0 \\ -1 \\ 1
	\end{pmatrix}
\end{equation*}
in the lattice $N = \ZZ^3$. They are the vertices of a lattice pentagon which lies on the horizontal plane $\RR^2 \times \{ 1\}$ in $N_\RR = \RR^3$.
Let $P$ be the convex hull of the $10$ points
\begin{gather*}
	\pm a, \ \pm b, \ \pm c, \ \pm d, \ \pm e
\end{gather*}
in $N$ (see Figure~\ref{fig:polytope_P}).
It is clear that $P$ is a centrally symmetric Fano polytope and has $10$ facets:
\begin{enumerate}
	\item[(i)] the pentagon $\conv{a,b,c,d,e}$ which lies on the horizontal plane $\RR^2 \times \{1\}$;
	\item[(i')] the pentagon $\conv{-a,-b,-c,-d,-e}$ which lies on the horizontal plane $\RR^2 \times \{-1\}$;
	\item[(ii)] the vertical rectangle $\conv{a,e,-c,-d}$;
		\item[(ii')] the vertical rectangle $\conv{-a,-e,c,d}$;
	\item[(iii)] the triangle $\conv{e,-b,-c}$;
	\item[(iii')] the triangle $\conv{-e,b,c}$;
	\item[(iv)] the triangle $\conv{a,b,-d}$;
	\item[(iv')] the triangle $\conv{-a,-b,d}$;
	\item[(v)] the triangle $\conv{-b,d,e}$;
	\item[(v')] the triangle $\conv{b,-d,-e}$.
\end{enumerate}
Moreover, one sees that the origin is the unique interior lattice point of $P$, so $P$ is a canonical polytope.
One can see that $P$ is $\mathrm{GL}_3(\ZZ)$-equivalent to the canonical $3$-dimensional polytope with ID 652083 in the Graded Ring Data Base \cite{kasprzyk_canonical}.

Consider the dual lattice $M = \Hom_\ZZ(N,\ZZ)$ and the polar polytope of $P$; this is the polytope $Q$ which is the convex hull of the points
\begin{gather*}
\pm (1,0,0), \qquad \quad
\pm (0,1,0), \\
\pm (0,0,1), \qquad \quad
\pm (1,-1,0), \\
\pm \! \left( \frac{2}{3}, \frac{2}{3}, - \frac{1}{3} \right)
\end{gather*}
 in $M_\RR$.
One can check that $Q$ has exactly $9$ lattice points, namely the origin and $\pm (1,0,0)$, 
$\pm (0,1,0)$, 
$\pm (0,0,1)$,
$\pm (1,-1,0)$.
The normalised volume of $Q$ is $\frac{40}{3}$.
The intersection of $Q$ with the horizontal plane $\RR^2 \times \{0\} \subset M_\RR$ is the hexagon with vertices
$\pm (1,0,0)$, $\pm(0,1,0)$, $\pm(-1,1,0)$ (see the green hexagon in Figure~\ref{fig:solo}).

 \section{The variety $X$} \label{sec:variety}

Let $X$ be the toric variety associated to the face fan (also called spanning fan) of $P$.
Since $P$ is a canonical $3$-dimensional polytope, $X$ is a Fano $3$-fold with canonical singularities.
The polytope $Q$ is the moment polytope associated to the toric boundary of $X$, which is an anticanonical divisor.
The anticanonical degree $(-K_X)^3$ of $X$ coincides with the normalised volume of $Q$, so it is $\frac{40}{3}$.
Since $Q$ is centrally symmetric, the origin is the barycentre of $Q$, thus $X$ is K-polystable by \cite{BJ,berman_polystability}.

The face fan of $P$ gives an affine open cover of $X$.
We denote by $U_{a,b,-d}$ the toric affine open subscheme of $X$ associated to the cone over the facet of $P$ with vertices $a,b,-d$, and similarly for the other facets.
Now we analyse the singularities which appear on $X$.
\begin{enumerate}
	\item[(i)] $U_{a,b,c,d,e}$ and $U_{-a,-b,-c,-d,-e}$ are isomorphic to the affine cone over the anticanonical embedding in $\PP^7$ of the smooth del Pezzo surface of degree $7$, i.e.\ the blowup of $\PP^2$ in $2$ points. This is an isolated Gorenstein singularity.
	\item[(ii)]  $U_{a,e,-c,-d}$ and $U_{-a,-e,c,d}$ are isomorphic to the hypersurface
	\[
	V = \Spec \CC[x,y,z,w] / (x^2 y^2 - zw),
	\]
	whose singular locus is $1$-dimensional and has $2$ irreducible components, which correspond to the edges with lattice length $2$.
	\item[(iii+iv)] $U_{e,-b,-c}$,  $U_{-e,b,c}$, 
$U_{a,b,-d}$ and
 $U_{-a,-b,d}$ are isomorphic to the hypersurface
 \[
 cA_1 = \Spec \CC[x,y,z,w] / (xy-z^2),
 \]
 whose singular locus is $1$-dimensional and irreducible and corresponds to the edge with lattice length $2$.
	
	\item[(v)] $U_{-b,d,e}$ and $U_{b,-d,-e}$ are isomorphic to the isolated cyclic quotient singularity $\frac{1}{3}(1,1,2)$. This is a $\QQ$-Gorenstein non-Gorenstein terminal singularity which is rigid, both with respect to flat deformations and to $\QQ$-Gorenstein deformations.
\end{enumerate}

The singular locus of $X$ has $6$ connected components:
\begin{enumerate}
	\item[(i)] the $0$-stratum of $U_{a,b,c,d,e}$;
	\item[(i')] the $0$-stratum of $U_{-a,-b,-c,-d,-e}$;
	\item[(ii+iii+iv)] the union of the singular loci of $U_{a,e,-c,-d}$, $U_{e,-b,-c}$, $U_{a,b,-d}$, which is the union of two smooth rational curves meeting transversally at the $0$-stratum of $U_{a,e,-c,-d}$;
		\item[(ii'+iii'+iv')] the union of the singular loci of $U_{-a,-e,c,d}$, $U_{-e,b,c}$, $U_{-a,-b,d}$, which is the union of two smooth rational curves meeting transversally at the $0$-stratum of $U_{-a,-e,c,d}$;
	\item[(v)] the $0$-stratum of $U_{-b,d,e}$;
		\item[(v')] the $0$-stratum of $U_{b,-d,-e}$.
\end{enumerate}

It is clear that the non-Gorenstein locus of $X$ consists of the $0$-strata of $U_{-b,d,e}$ and $U_{b,-d,-e}$.
Let $\epsi \colon \frakX \to X$ be the index $1$ cover stack of $X$.
$\frakX$ is a Deligne--Mumford stack of finite type over $\CC$, $\frakX$ is Gorenstein and $\epsi$ is an isomorphism over the Gorenstein locus of $X$.
Therefore $\epsi$ is an isomorphism away from the $0$-strata of $U_{-b,d,e}$ and $U_{b,-d,-e}$.
One immediately sees that $\frakX$ is lci away from the $0$-strata of $U_{a,b,c,d,e}$ and $U_{-a,-b,-c,-d,-e}$.

\section{Some computations} \label{sec:computation}

\subsection{Cotangent sheaves and groups}

For each $i=0,1,2$, consider the coherent sheaves $\cT^i_X = \cExt^i_X(\Omega_X, \cO_X)$ and $\cTqG{i}{X} = \epsi_\star \cExt^i_X(\Omega_\frakX, \cO_X)$ on $X$, and the $\CC$-vector spaces $\TT^i_X = \Ext^i_X(\Omega_X, \cO_X)$ and $\TTqG{i}{X} = \Ext^i_\frakX(\Omega_\frakX, \cO_\frakX)$. All of them are $M$-graded, because the torus $T = N \otimes_\ZZ \GG_\mathrm{m} = \Spec \CC[M]$ acts on $X$.

The sheaves $\cT^0_X$ and $\cTqG{0}{X}$ coincide everywhere on $X$. For $i>0$ the sheaves $\cT^i_X$ and $\cTqG{i}{X}$ coincide on the Gorenstein locus of $X$; therefore when considering restrictions to Gorenstein open subschemes we will suppress the superscript $\mathrm{qG}$. For $i>0$, the set-theoretical support of $\cTqG{i}{X}$ is contained in the singular locus of $X$.

\subsection{The sheaf $\cTqG{1}{X}$} \label{sec:T1}
We analyse the restriction of $\cTqG{1}{X}$ to the affine charts given by the facets of $P$.
\begin{enumerate}
	\item[(i)] The restriction to $U_{a,b,c,d,e}$ is supported on the $0$-stratum of $U_{a,b,c,d,e}$ and its global sections are $\CC^2$ in degree $(0,0,-1) \in M$ by \cite{altmann_computation_tangent}.

	\begin{figure}
	\centering
	\includegraphics[scale = 0.45]{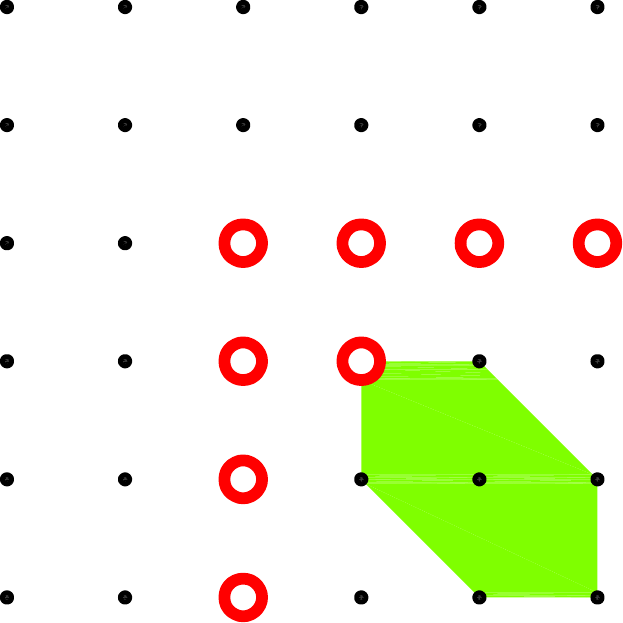}
	\caption{A portion of the horizontal plane $\RR^2 \times \{ 0\}$ in $M_\RR = \RR^3$. 
		The green hexagon is the slice of the polytope $Q$.
		The centre of the hexagon is the origin.
		The degrees of the homogeneous components of $\TT^1_{U_{a,e,-c,-d}}$ are denoted by red circles.}
	\label{fig:solo}
\end{figure}

	\item[(ii)] If one considers the hypersurface $V = \Spec \CC[x,y,z,w] / (x^2 y^2 - zw)$, one has that $\TT^1_V$ is the $\cO_V$-module
	\[
	\CC[x,y,z,w] / (xy^2, x^2 y, z, w) = \CC 1 \oplus \CC zw \oplus \bigoplus_{n \geq 1} \left( \CC z^n \oplus \CC w^n \right).
	\]
	The homogeneous components of $\TT^1_{U_{a,e,-c,-d}}$ with respect to the $M$-grading are:
	\begin{equation*}
	\TT^1_{U_{a,e,-c,-d}}(m) = \begin{cases}
			\CC \quad &\text{if } m = (-1,1,0) \text{ or }m = (-2,2,0), \\
	\CC \quad &\text{if } m = (-2+n,2,0) \text{ and } n \geq 1, \\
		\CC \quad &\text{if } m = (-2,2-n,0) \text{ and } n \geq 1, \\ 
		0 	 \quad &\text{otherwise.} 
	\end{cases}
	\end{equation*}
See Figure~\ref{fig:solo}.

	\item[(iii)] If one considers the hypersurface $cA_1 = \Spec \CC[x,y,z,w] / (xy - z^2)$, one has that $\TT^1_{cA_1}$ is the $\cO_{cA_1}$-module
	\[
	\CC[x,y,z,w] / (x,y,z) = \bigoplus_{n \geq 0} \CC w^n.
	\]
		The homogeneous components of $\TT^1_{U_{e,-b,-c}}$ with respect to the $M$-grading are:
	\begin{equation*}
		\TT^1_{U_{e,-b,-c}}(m) = \begin{cases}
			\CC \quad &\text{if } m = (-2-n,2,0) \text{ and } n \geq 0, \\
			0 	 \quad &\text{otherwise.} 
		\end{cases}
	\end{equation*}
	See the left of Figure~\ref{fig:accoppiati_0cocycle}.

\begin{figure}
	\centering
	\includegraphics[scale = 0.45]{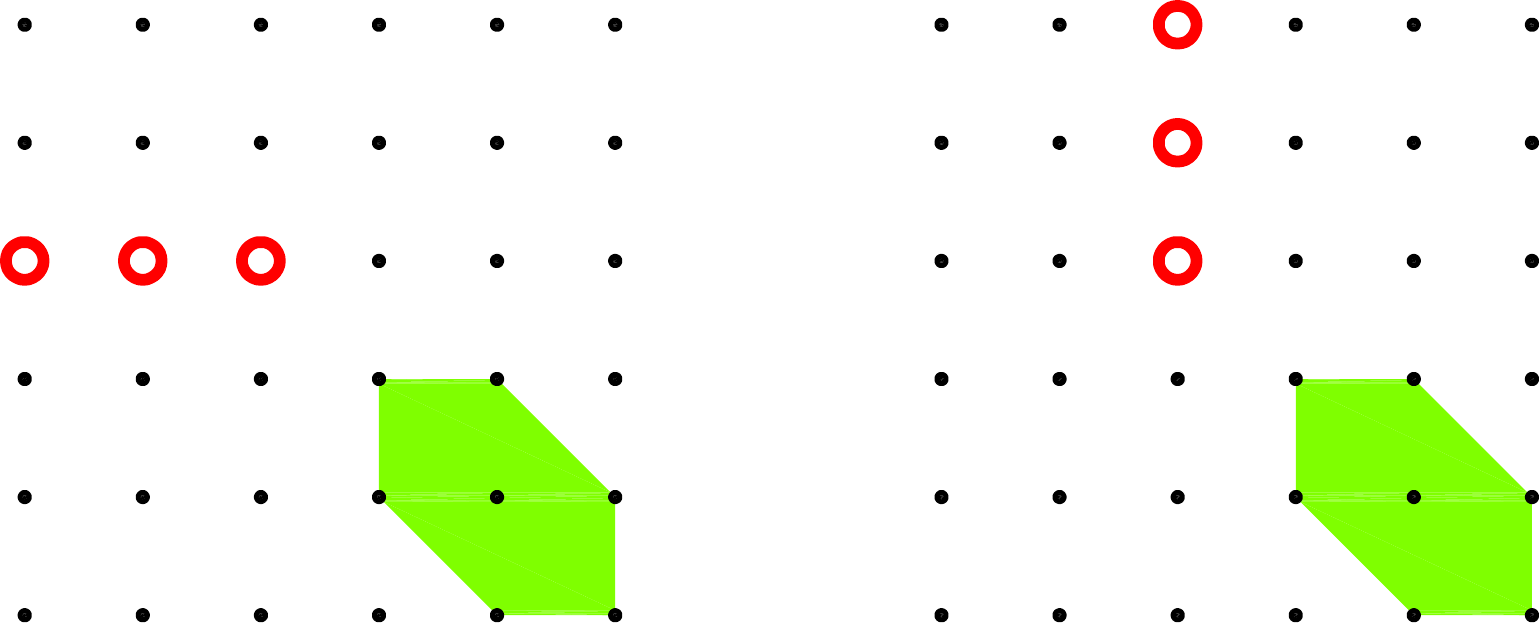}
	\caption{The red circles denote the degrees 
		of the homogeneous components of $\TT^1_{U_{e,-b,-c}}$ (left)  and of $\TT^1_{U_{a,b,-d}}$ (right).}
	\label{fig:accoppiati_0cocycle}
\end{figure}

	\item[(iv)] 
			The homogeneous components of $\TT^1_{U_{a,b,-d}}$ with respect to the $M$-grading are:
	\begin{equation*}
		\TT^1_{U_{a,b,-d}}(m) = \begin{cases}
			\CC \quad &\text{if } m = (-2,2+n,0) \text{ and } n \geq 0, \\
			0 	 \quad &\text{otherwise.} 
		\end{cases}
	\end{equation*}
	See the right of Figure~\ref{fig:accoppiati_0cocycle}.
	
	\item[(v)] The restriction to $U_{-b,d,e}$ is zero.

\end{enumerate}
The discussion about (i'), (ii'), (iii'), (iv'), (v') is completely analogous and is omitted: it is enough to apply the reflection $-\mathrm{id}$.

Now we want to understand the restriction of the sheaf $\cTqG{1}{X}$ to the double intersections of the toric charts given by the facets of $P$. We will ignore the smooth ones, because there the sheaf $\cTqG{1}{X}$ vanishes.
\begin{enumerate}
	\item[(ii$\cap$iii)] The intersection $U_{a,e,-c,-d} \cap U_{e,-b,-c}$ is isomorphic to the hypersurface $\Spec \CC[x,y,z,w^\pm] / (xy-z^2)$, so its $\TT^1$-module is isomorphic to
	\[
	\CC[x,y,z,w^\pm] / (x,y,z) = \bigoplus_{n \in \ZZ} \CC w^n.
		\]
		The homogeneous components with respect to the $M$-grading are:
		\begin{equation*}
			\TT^1_{U_{a,e,-c,-d} \cap U_{e,-b,-c}}(m) = \begin{cases}
				\CC \quad &\text{if } m = (-2+n,2,0) \text{ and } n \in \ZZ, \\
				0 	 \quad &\text{otherwise.} 
			\end{cases}
		\end{equation*}
	See the left of Figure~\ref{fig:accoppiati_double_intersections}.

	\begin{figure}
	\centering
	\includegraphics[scale = 0.45]{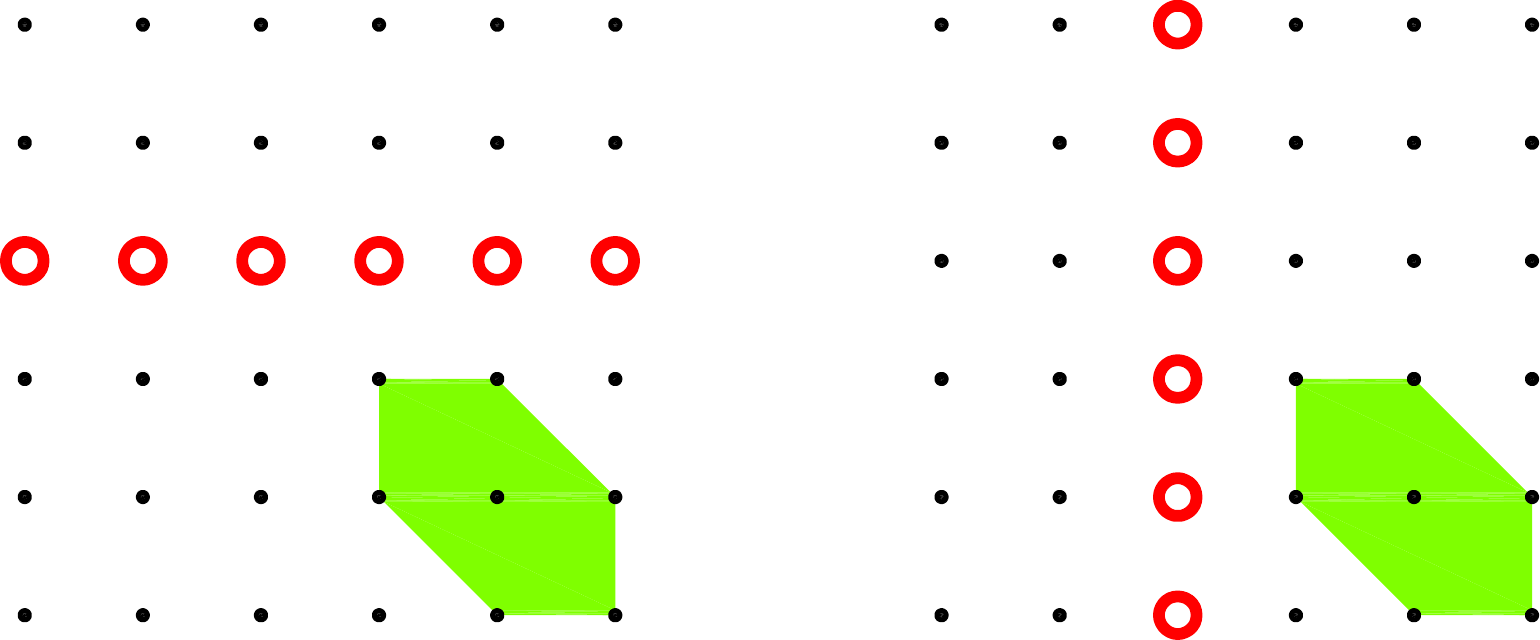}
	\caption{The red circles denote the degrees 
		of the homogeneous components of $\TT^1_{U_{a,e,-c,-d} \cap U_{e,-b,-c}}$ (left)  and of $\TT^1_{U_{a,e,-c,-d} \cap U_{a,b,-d}}$ (right).}
	\label{fig:accoppiati_double_intersections}
\end{figure}
	
	\item[(ii$\cap$iv)] 
	The intersection $U_{a,e,-c,-d} \cap  U_{a,b,-d}$ is isomorphic to $\Spec \CC[x,y,z,w^\pm] / (xy-z^2)$.
	The $M$-homogeneous components of the restriction of $\cTqG{1}{X}$ to the intersection $U_{a,e,-c,-d} \cap U_{a,b,-d}$ are
			\begin{equation*}
		\TT^1_{U_{a,e,-c,-d} \cap  U_{a,b,-d}}(m) = \begin{cases}
			\CC \quad &\text{if } m = (-2,2+n,0) \text{ and } n \in \ZZ, \\
			0 	 \quad &\text{otherwise.} 
		\end{cases}
	\end{equation*}
	See the right of Figure~\ref{fig:accoppiati_double_intersections}.
\end{enumerate}
The discussion about the intersections (ii'$\cap$iii') and (ii'$\cap$iv') is omitted because it can easily obtained from above by applying the reflection $-\mathrm{id}$.

Now we have all information to construct the \v{C}ech complex of the sheaf $\cTqG{1}{X}$ with respect to the affine cover given by the charts associated to the facets of $P$. Notice that, since triple intersections are smooth, this complex is concentrated in degrees $0$ and $1$:
\[
C^0 \overset{d}{\longrightarrow} C^1.
\]
We analyse its homogeneous components with respect to the $M$-grading.
\begin{enumerate}[label=$\bullet$]
	\item If $m=\pm(0,0, 1)$, then $d(m)$ is the zero map $\CC^2 \to 0$.
	\item If $m=\pm(-1,1,0)$, then $d(m)$ is the zero map $\CC \to 0$.
	\item If $m=\pm(-2,2,0)$, then $d(m)$ is
	\[
	\CC^3 \xrightarrow{\begin{pmatrix}
		1 & -1 & 0 \\ 1 & 0 & -1
		\end{pmatrix}}
	\CC^2.
	\]
	\item If $m$ is none of the degrees above, then $d(m)$ is the identity of $0$, of $\CC$ or of $\CC^2$.
\end{enumerate}
We have proved that $d$ is surjective, therefore $\rH^1(X, \cTqG{1}{X}) = 0$.
Furthermore, we have computed the homogeneous components of $\rH^0(X, \cTqG{1}{X})$:
			\begin{equation*}
\rH^0(X, \cTqG{1}{X})(m) = \begin{cases}
		\CC^2 \quad &\text{if } m = \pm (0,0,1), \\
		\CC \quad &\text{if } m = \pm (-1,1,0), \\
		\CC \quad &\text{if } m = \pm (-2,2,0), \\
		0 	 \quad &\text{otherwise.} 
	\end{cases}
\end{equation*}
In particular $\dim \rH^0(X, \cTqG{1}{X}) = 8$.

\subsection{The hull of $\DefqG{X}$} \label{sec:hull}
Let $\DefqG{X}$ denote the functor of $\QQ$-Gorenstein deformations of $X$.
Let $A$ denote its hull.
In other words, the base of the miniversal $\QQ$-Gorenstein deformation of $X$ is the formal spectrum $\Spf A$.

Since $\rH^1(X, \cT^0_X) = 0$ and $\rH^2(X, \cT^0_X) = 0$ by \cite[Proof of Theorem~5.1]{totaro} (see also \cite[\S14.4.3]{petracci_survey}), the natural homomorphism $\TTqG{1}{X} \to \rH^0(X, \cTqG{1}{X})$ is an isomorphism.
From the conclusion of \S\ref{sec:T1} we deduce that $\TTqG{1}{X}$ has dimension $8$ and we can choose coordinates $t_1, \dots, t_8$ with weights $(0,0,1)$, $(0,0,1)$, $(0,0,-1)$, $(0,0,-1)$, $(-1,1,0)$, $(-2,2,0)$, $(1,-1,0)$, $(2,-2,0)$ in $M$, respectively.
This implies that $A$ is a quotient of the power series ring $\CC \pow{t_1, \dots, t_8}$ with respect to an ideal which is contained in $(t_1, \dots, t_8)^2$.

In order to understand this ideal we need to compute the obstructions.
Since $\frakX$ is lci away from the $0$-strata of $U_{a,b,c,d,e}$ and $U_{-a,-b,-c,-d,-e}$,
the sheaf $\cTqG{2}{X}$ is set-theoretically supported on these two points, hence
\begin{equation*}
\rH^0(X, \cTqG{2}{X}) = \TT^2_{U_{a,b,c,d,e}} \oplus \TT^2_{U_{-a,-b,-c,-d,-e}}.
\end{equation*}
Moreover, the vanishing of $\rH^2(X, \cT^0_X)$ and of $\rH^1(X, \cTqG{1}{X})$ implies that the natural  homomorphism $\TTqG{2}{X} \to \rH^0(X, \cTqG{2}{X})$ is an isomorphism.

We deduce that the product of restriction maps
\begin{equation} \label{eq:map_restriction}
\DefqG{X} \longrightarrow \Def{U_{a,b,c,d,e}} \times \Def{U_{-a,-b,-c,-d,-e}}.
\end{equation}
induces and isomorphism on tangent spaces and on obstruction spaces. Therefore this map is smooth.

By \cite{altmann_versal} $\CC \pow{u,v} / (u^2, uv)$ is the hull of the deformation functor of the affine cone over the anticanonical embedding in $\PP^7$ of the smooth del Pezzo surface of degree $7$.
From the smoothness of   \eqref{eq:map_restriction} we deduce
\begin{equation} \label{eq:A}
A = \CC \pow{t_1, \dots, t_8} / (t_1^2, t_1 t_2, t_3^2, t_3 t_4).
\end{equation}
This shows that the base of the miniversal $\QQ$-Gorenstein deformation of $X$ is irreducible and non-reduced.

\subsection{The action the automorphism group}

The automorphism group of the polytope $P$ is generated by the two involutions
\begin{equation*}
\sigma = \begin{pmatrix}
	-1 & 0 & 0 \\
	0 & -1 & 0 \\
	0 & 0 & -1
\end{pmatrix}
\quad \text{ and } \quad
\tau = \begin{pmatrix}
	0 & 1 & 0 \\
	1 & 0 & 0 \\
	0 & 0 & 1
\end{pmatrix}
\end{equation*}
which commute in $\GL(N) = \GL_3(\ZZ)$.
Let $G$ denote the automorphism group of $X$ and let $T = N \otimes_\ZZ \GG_\mathrm{m} = \Spec \CC[M]$ be the torus acting on $X$.
Since every facet of $Q$ has no interior lattice point, 
by \cite[Proposition~2.8]{ask_petracci} $G$ is the semidirect product $T \rtimes (C_2 \times C_2)$, where $C_2$ denotes the group of order $2$.
The algebraic group $G$ acts formally on $A$ as follows:
\begin{itemize}
	\item the torus $T$ formally acts linearly on $A$ via the weights of the $T$-representation $\TTqG{1}{X}$, hence $t_1, \dots, t_8$ have weights $(0,0,1)$, $(0,0,1)$, $(0,0,-1)$, $(0,0,-1)$, $(-1,1,0)$, $(-2,2,0)$, $(1,-1,0)$, $(2,-2,0)$, respectively;
	\item $\sigma$ acts by permuting the coordinates $t_1, \dots, t_8$ as follows: $t_1 \leftrightarrow t_3$, $t_2 \leftrightarrow t_4$, $t_5 \leftrightarrow t_7$, $t_6 \leftrightarrow t_8$;
	\item $\tau$ permutes the coordinates $t_1, \dots, t_8$ as follows: $t_1$, $t_2$, $t_3$, $t_4$ are left fixed and $t_5 \leftrightarrow t_7$, $t_6 \leftrightarrow t_8$.
\end{itemize}
We need to compute the invariant subring $A^G$, which is the limit of $(A / \frakm_A^{i+1})^G$ as $i \geq 0$, where $\frakm_A$ denotes the maximal ideal of $A$.
Actually there is no harm in considering polynomials instead of power series
\[
B = \CC[t_1, \dots, t_8] / (t_1^2, t_1 t_2, t_3^2, t_3 t_4)
\]
and on $B$ there is an honest, i.e.\ non-formal, action of the algebraic group $G$.
In the $G$-action the two sets of variables $\{t_1, t_2, t_3, t_4\}$ and $\{t_5, t_6, t_7, t_8\}$ are not mixed. So we can consider
\begin{equation*}
R = \CC[t_1, t_2, t_3, t_4] / (t_1^2, t_1 t_2, t_3^2, t_3 t_4) \quad \text{ and } \quad S = \CC[t_5, t_6, t_7, t_8]
\end{equation*}
and study the $G$-actions on $R$ and on $S$ separately.

The invariant subring $\CC[t_1, t_2, t_3, t_4]^T$ is generated by $x_1 = t_1 t_3$, $x_2 = t_2 t_4$, $x_3 = t_2 t_3$ and $x_4 = t_1 t_4$, which satisfy the relation $x_1 x_2 - x_3 x_4  = 0$.
\[
\xymatrix{
& \CC[t_1, t_2, t_3, t_4] \ar@{->>}[r] & R \\
\CC[x_1, x_2, x_3, x_4] \ar@{->>}[r] & \CC[t_1, t_2, t_3, t_4]^T \ar@{->>}[r] \ar@{_{(}->}[u] & R^T \ar@{_{(}->}[u]
}
\]
One can show that the kernel of $\CC[x_1, x_2, x_3, x_4] \onto R^T$ is generated by
$x_1^2 = t_1^2 t_3^2$, $x_3^2 = t_2^2 t_3^2$, $x_4^2 = t_1^2 t_4^2$, $x_1 x_2 = x_3 x_4 = t_1 t_2 t_3 t_4$, $x_1 x_3 = t_1 t_2 t_3^2$, $x_1 x_4 = t_1^2 t_3 t_4$, $x_2 x_3 = t_2^2 t_3 t_4$, $x_2 x_4 = t_1 t_2 t_4^2$,
i.e.\ by all degree $2$ monomials in $x_1, x_2, x_3, x_4$ with the exception of $x_2^2$.
Therefore we have an isomorphism
\[
R^T \simeq \CC[x_1, x_2, x_3, x_4] / (x_1^2, x_3^2, x_4^2, x_1 x_2, x_1 x_3, x_1 x_4, x_2 x_3, x_2 x_4, x_3 x_4).
\]
The involution $\sigma$ keeps $x_1$ and $x_2$ fixed and swaps $x_3$ and $x_4$.
The involution $\tau$ acts trivially on $R$, and therefore also on $R^T$.
Hence $R^G = (R^T)^\sigma$ is generated by $y_1 = x_1$, $y_2 = x_2$, $y_3 = x_3 + x_4$, $y_4 = x_3 x_4$.
\[
\xymatrix{
	& \CC[x_1, x_2, x_3, x_4] \ar@{->>}[r] & R^T \\
	\CC[y_1, y_2, y_3, y_4] \ar@{=}[r] & \CC[x_1, x_2, x_3, x_4]^\sigma \ar@{->>}[r] \ar@{_{(}->}[u] & R^G \ar@{_{(}->}[u]
}
\]
The kernel of $\CC[y_1, y_2, y_3, y_4] \onto R^G$ is generated by $y_1^2$, $y_1 y_2$, $y_4$, $y_3^2$, $y_1 y_3$, $y_2 y_3$.
Therefore we have isomorphisms
\begin{equation} \label{eq:R^G}
R^G \simeq \frac{\CC[y_1, y_2, y_3, y_4]}{(y_1^2, y_1 y_2, y_4, y_3^2, y_1 y_3, y_2 y_3)} \simeq \frac{\CC[y_1, y_2, y_3]}{(y_1^2, y_3^2, y_1 y_2, y_1 y_3, y_2 y_3)}.
\end{equation}

The invariant subring $S^T$ is generated by $z_1 = t_5^2 t_6$, $z_2 = t_5 t_7$, $z_3 = t_7^2 t_8$, $z_4 = t_6 t_8$, which satisfy the relation $z_1 z_3 - z_2^2 z_4 = 0$.
Therefore we have an isomorphism
\[
S^T \simeq \CC[z_1, z_2, z_3, z_4] / (z_1 z_3 - z_2^2 z_4).
\]
The action of $\tau$ on $S$ coincides with the action of $\sigma$ on $S$, therefore we can ignore $\tau$.
One immediately sees that $\sigma$ swaps $z_1$ and $z_3$ and leaves $z_2$ and $z_4$ fixed.
This implies that $S^G = (S^T)^\sigma$ is generated by $v_1 = z_1 + z_3$, $v_2 = z_1 z_3$, $v_3 = z_2$, $v_4 = z_4$.
\[
\xymatrix{
	& \CC[z_1, z_2, z_3, z_4] \ar@{->>}[r] & S^T \\
	\CC[v_1, v_2, v_3, v_4] \ar@{=}[r] & \CC[z_1, z_2, z_3, z_4]^\sigma \ar@{->>}[r] \ar@{_{(}->}[u] & S^G \ar@{_{(}->}[u]
}
\]
The kernel of $\CC[v_1, v_2, v_3, v_4] \onto S^G$ is generated by $v_2 - v_3^2 v_4$.
Therefore we have isomorphisms
\begin{equation} \label{eq:S^G}
S^G \simeq \frac{\CC[v_1, v_2, v_3, v_4]}{(v_2 - v_3^2 v_4)} \simeq \CC[v_1, v_3, v_4].
\end{equation}

By \eqref{eq:R^G} and \eqref{eq:S^G}, by taking the completion of $B^G \simeq R^G \otimes_\CC S^G$, and by changing the names of the variables, we get an isomorphism
\begin{equation} \label{eq:A^G}
A^G \simeq \CC \pow{u_1, \dots, u_6} / (u_1^2,  u_2^2, u_1 u_2, u_1 u_3, u_2 u_3).
\end{equation}

\subsection{Conclusion}

By the Luna \'etale slice theorem for algebraic stacks \cite{luna_etale_slice_stacks} there exists a commutative diagram
\[
\xymatrix{
[\Spf A \ / G] \ar[d] \ar[r] & \stack{3} \ar[d] \\
\Spf A^G \ar[r] & \modspace{3}
}
\]
where the horizontal arrows are formally \'etale and map the closed point to the point corresponding to $X$.
The isomorphisms \eqref{eq:A} and \eqref{eq:A^G}
imply Theorem~\ref{thm}.


\section{Mirror symmetry}
\label{sec:mirror_symmetry}

Let us consider Laurent polynomials in $3$ variables $x$, $y$, $z$ whose Newton polytope is $P$.
For such $f \in \CC[x^\pm, y^\pm, z^\pm]$ we insist that
\begin{itemize}
\item the coefficient of the monomial $1$ (which corresponds to the origin) is $0$,
\item the coefficients of the monomials corresponding to vertices of $P$ are equal to $1$,
\end{itemize}
then we still need to decide the coefficients of the monomials $z$, $z^{-1}$, $x$, $x^{-1}$, $y$, $y^{-1}$.
By imposing that the restriction of $f$ to each facet of $P$ is $0$-mutable in the sense of \cite{cofipe} (see also \cite{maximally_mutable}), we get the Laurent polynomial
\begin{align*}
	f &= (x+y+xy)(1 + x^{-1} y^{-1})z +
	(x^{-1} + y^{-1} + x^{-1}y^{-1})(1 + xy) z^{-1} + \\
&+ 2(x+x^{-1} + y + y^{-1})
\end{align*}
which can also be written as
\begin{align*}
f &= (x+x^{-1} + y + y^{-1})z^{-1} (1+z)^2 +z(1+xy) + z^{-1} (1+x^{-1} y^{-1}).
\end{align*}

One can consider the  \emph{classical period} of $f$, which is a power series $\pi_f \in \CC \pow{t}$ defined in terms of variations of Hodge structures of the fibration $f \colon \GG_\mathrm{m}^3 \to \AA^1$. We refer the reader to \cite{galkin_usnich, mirror_symmetry_and_fano_manifolds} for the precise definition of $\pi_f$. In this particular case the first terms of $\pi_f$ are:
\begin{equation*}
	\pi_f = 1+ 28 t^2 + 144 t^3 + 3516 t^4 + 38400 t^5 + \cdots.
\end{equation*}

From the results in \S\ref{sec:hull} it is easy to prove that the general fibre of the miniversal $\QQ$-Gorenstein deformation of $X$ is a terminal Fano $3$-fold $X'$ with two $\frac{1}{3}(1,1,2)$ points.
According to \cite{mirror_symmetry_and_fano_manifolds} (see also \cite{chp,log_crepant1, petracci_roma}), it is conjectured that $X'$ is \emph{mirror} to $f$, in the sense that $\pi_f$ coincides with \emph{regularised quantum period} of $X'$. Here the regularised quantum period $\widehat{G}_{X'}$ is a certain generating function for genus $0$ Gromov--Witten invariants of $X'$; we refer the reader to \cite[Definition~3.2]{oneto_petracci} for the precise definition.

Unfortunately we are not able to prove the equality $\widehat{G}_{X'} = \pi_f$ because we do not know how to compute Gromov--Witten invariants of $X'$; indeed, we do not know a presentation of $X'$ as a complete intersection in a toric variety or in a flag variety.
Nonetheless the reader should appreciate that this type of mirror symmetry conjectures establishes an intriguing correspondence between algebraic geometry (i.e.\ the connected component of moduli containing $X$ and $X'$, and the Gromov--Witten theory of $X'$) and combinatorics (i.e.\ the polynomial $f$ and its classical period $\pi_f$).

\bibliography{Biblio_agplus}

\end{document}